\documentclass{elsarticlenonels}

\usepackage{amssymb}
\usepackage{amsmath}
\usepackage{amsfonts}
\usepackage[dvipsnames,usenames]{color}
\usepackage[colorlinks=true,
pdfstartview=FitV, linkcolor=cyan,
citecolor=cyan, urlcolor=cyan]{hyperref}
\newtheorem{remark}{Remark}
\newtheorem{proposition}{Proposition}
\newtheorem{theorem}{Theorem}
\newtheorem{lemma}{Lemma}
\newtheorem{definition}{Definition}
\numberwithin{theorem}{section}
\numberwithin{proposition}{section}
\numberwithin{lemma}{section}
\numberwithin{remark}{section}

\catcode`\=\active \def {\`e}        
\catcode`\Ž=\active \def Ž{\'e}        
\catcode`\ˆ=\active \def ˆ{\`a}        
\catcode`\"=\active \def "{\`{\i}}        
\catcode`\˜=\active \def ˜{\`o}        
\catcode`\=\active \def {\`u}        
\catcode`\"=\active \def "{\^{\i}}       

\def\persone#1{\textsc{\textcolor{BrickRed}{#1}}}

\def\spray{{\BS}}

\def\torstensor{\persone{tors}}

\def\brack#1#2{[#1,#2]}

\def\endprova{\hfill{$\blacksquare$}\goodbreak\medskip}

\def\BUNDLE{\BBB}

\def\Null#1{\mathbf{ker}\di{#1}}
\def\Range#1{\mathbf{im}\di{#1}}
\def\Bproj{\Bp}
\def\Tproj{\Btau}

\def\proof{\goodbreak\par\noindent\textbf{Proof.} }

\def\duetonde#1#2{{( \! (}\kern.1em
#1,#2\kern.1em{) \! )}}

\def\TE{{\TANG\EEE}}
\def\VE{\VERT\EEE}
\def\HE{\HORI\EEE}

\def\Tder#1{\TT_{#1}}

\def\VPROJ{\mathrm{P_V}}
\def\HPROJ{\mathrm{P_H}}

\def\Hlift{\BH}

\def\FL#1#2{\mathbf{Fl}^{#1}_{#2}}

\def\TprojM{\Tproj_\MMM}
\def\TprojN{\Tproj_\NNN}
\def\TprojE{\Tproj_\EEE}

\def\path{\Bgamma}

\def\forw{{\Uparrow\,}}

\def\push{{\uparrow}}
\def\pull{{\downarrow}}
\def\pttt{\lambda}
\def\ptau{\mu}

\def\Lieder{\cL}

\def\BBB{\mathbb{B}}
\def\map{{\Bvarphi}}

\def\Bvtau{{\Bv_{\tau}}}

\def\curvoperator{\textsc{curv}}
\def\curvopbar{\segno{\textsc{curv}}}

\def\curvform{\mathbf{R}}
\def\cocurvform{\mathbf{R}^{\mathbf{c}}}

\def\parder#1#2{\partial_{#1=#2}\,}

\let\segno=\overline

\def\EEE{\mathbb{E}}
\def\FFF{\mathbb{F}}

\def\TANG{\mathbb{T}}

\def\VERT#1{\mathbb{V}#1}
\def\HORI#1{\mathbb{H}#1}

\def\MMM{\mathbb{M}}

\def\NNN{\mathbb{N}}

\def\ID#1{\mathbf{id}_{\,#1}}

\def\TM{{\TANG\MMM}}

\def\TTANG{\TANG\TANG}
\def\TTM{\TTANG\MMM}

\def\TN{\TANG\NNN}

\def\nablabar{\overline\nabla}

\def\inv#1{#1^{-1}}

\let\aa=\alpha
\let\bb=\beta

\def\cont{\mathrm{C}}

\def\di#1{(#1)}

\def\coppia#1#2{(#1\,,#2)}

\def\BSigma{\boldsymbol{\Sigma}}
\def\Btau{\boldsymbol{\tau}}

\def\Bpi{\boldsymbol{\pi}}

\def\Bss{\boldsymbol{\sigma}}

\def\Bgamma{\boldsymbol{\gamma}}

\def\Bvarphi{\boldsymbol{\varphi}}

\def\torstensor{\persone{tors}}

\def\Bb{\mathbf{b}}

\def\Be{\mathbf{e}}

\def\Bp{\mathbf{p}}

\def\Bs{\mathbf{s}}

\def\Bu{\mathbf{u}}
\def\Bv{\mathbf{v}}

\def\Bx{\mathbf{x}}

\def\BS{\mathbf{S}}

\def\BX{\mathbf{X}}
\def\BY{\mathbf{Y}}

\def\BH{\mathbf{H}}

\def\sub#1{{}_{\lower2pt\hbox{$\scriptstyle#1$}}}

\def\cL{\mathcal{L}}

\def\cF{\mathcal{F}}

\def\II{I}

\def\TT{T}

\def\ff{f}

\def\cc{c}
\def\ttt{t}

\def\equaldef{:=}
\def\perogni{\quad\forall\,}
\def\Linmap#1{\textrm{\textit{BL}}\,\di{#1}}
\def\sp{\,;}

\def\equi{\,\Longleftrightarrow\,}
\def\implies{\,\Longrightarrow\,}
\def\scalar#1#2{{\langle}\kern.1em#1,#2\kern.1em{\rangle}}

\let\punto=\cdot

\def\set#1{\{#1\}}

\def\barra#1#2{ \bigg\vert _{\hbox{$\scriptstyle#1$}}^{\raise8.2pt\hbox{$\scriptstyle#2$}}  \, }

\def\lineare#1{[#1]}

\def\duetonde#1#2{{( \! (}\kern.1em\,#1\,,\,#2\,\kern.1em{) \! )}}

\def\with#1{^{\raise.5pt\hbox{$\scriptstyle#1$}}}

\def\duescalar#1#2{{\langle\langle}\kern.1em#1,#2\kern.1em{\rangle\rangle}}

\def\Ehresmann{\persone{Ehresmann}}
\def\Frobenius{\persone{Frobenius}}
\def\Hodge{\persone{Hodge}}
\def\Lie{\persone{Lie}}
\def\Leibniz{\persone{Leibniz}}
\def\Jacobi{\persone{Jacobi}}

\begin{document}
\begin{frontmatter}



\title{\persone{On natural derivatives and \\ the curvature formula in fibre bundles}}


\author[label1]{\persone{Giovanni Romano}}

\address[label1]{
University of Naples Federico II,\\
Department of Structural Engineering, \\
via Claudio 21, 80125 - Naples, Italy\\
e-mail: romano@unina.it}

\begin{abstract}
In a fibre bundle, natural derivatives of a section are defined as 
tangent vector fields on the image of a section of the fibre bundle.
A local extension to vector fields in the tangent bundle
leads to a direct proof of the formula expressing the curvature of
a connection in terms of covariant derivatives.
The result is based on a tensoriality argument
and extends to nonlinear connections on fibre bundles 
a well-known formula for linear connections on vector bundles.
\end{abstract}

\begin{keyword}
fibre bundles, \sep natural derivatives  \sep connections \sep integrability
\sep covariant derivatives.
\end{keyword}
\end{frontmatter}


\section{Introduction}

The notion of connection on a fibre bundle
was introduced by
\persone{Charles} \Ehresmann\ \cite{Ehresmann} in $\,1950\,$
and investigated by 
\persone{Paulette} \persone{Libermann} in
\cite{Libermann1969},
\cite{Libermann1973}, \cite{Libermann1982}.
Standard references on this topic are the article
by \persone{Kobayashi} 
\cite{Kobayashi1957} and the text 
\cite{KobayashiNomizu1963}.
The analysis developed in this paper makes also
reference to the treatment of the matter presented 
in \cite{KMS1993}.
Let us recall some well-known facts.
In the tangent bundle to a fibre bundle
the vertical distribution
is naturally defined by considering, at each point
of the manifold, the vectors tangent to the fibre through that point.
The vertical distribution is always integrable and the leaves
of the induced foliation are the fibres themselves.
The general definition of a connection as a (regular)
field of projections on the vertical subspaces of the tangent spaces
to a fibre bundle,
splits each tangent space into two complementary subspaces,
the vertical and the horizontal ones.
This leads naturally to the question about 
integrability of the horizontal distribution.
The involutivity condition provided by \Frobenius\ theorem
leads to the definition of the curvature
as obstruction against integrability 
of the horizontal distribution \cite{KMS1993}.
In this context I provide a new result, stated hereafter in Theorem \ref{th: curvnice}.
This result builds a direct bridge between 
the expression of the curvature in terms of horizontal lifts,
which is the one naturally steaming out of \Frobenius\ involutivity
condition, and the expression of the curvature in terms of covariant
derivatives, more suitable for applications.
The expression, which is well-known for linear connections on vector or principal bundles,
is extended by the new result to
general connections on fibre bundles.
The proof is based on the novel definition of natural derivative vector fields,
on an extension to a vector field in the tangent bundle
and on a direct, powerful tensoriality argument.
The analysis moves along the same line of thought
as for instance the one declared
in \cite{MangiarottiModugno1984}, by trying to avoid
unnecessary recourse to additional geometric structures.
In this respect the assumptions and the result of our 
Theorem \ref{th: curvnice} should be compared with the ones in 
\cite{KobayashiNomizu1963} Chapter III Theorem 5.1,
in \cite{ChoquetDeWittDillard1989} Chapter V-bis Section A.5, 
in \cite{MangiaSarda2000} Chapter 2 Section 2.4,
and in \cite{Michor2007} Corollary $\,19.16\,$, dealing with the
curvature of linear connections on vector bundles.

\section{Connection on a fibre bundle}
\label{sec: connection}
\index{\Ehresmann\ connection}

Let us recall some definitions and notations 
\cite{Saunders1989}, \cite{Lang1995}, \cite{RomanoDiff2007}.
Given two differentiable manifolds $\,\MMM\,$,$\,\NNN\,$,
the related tangent bundles with projections
$\,\TprojM\in\cont^1\di{\TM\sp\MMM}\,$,
$\,\TprojN\in\cont^1\di{\TN\sp\NNN}\,$
and a morphism $\,\map\in\cont^1\di{\MMM\sp\NNN}\,$,
a vector field $\,\BX\in\cont^1\di{\map\di\MMM\sp\TN}\,$
is $\,\map$-related to 
a vector field $\,\Bv\in\cont^1\di{\MMM\sp\TM}\,$ if 
$\,\BX\circ\map=\TT\map\circ\Bv\,$ where $\,\TT\,$ is the tangent functor.
For a diffeomorphism $\,\map\in\cont^1\di{\MMM\sp\NNN}\,$
the push and pull operations are then defined by
$\,\BX=\map\push\Bv\,$ and $\,\Bv=\map\pull\BX\,$.
The usual notation is 
$\,\map\push=\map_*\,$ and $\, \map\pull=\map^*\,$
but then too many stars do appear in the geometrical sky 
(push, duality, \Hodge\ star).
A fibre bundle is a surjective submersion
$\,\Bp\in\cont^1\di{\EEE\sp\MMM}\,$ with $\,\EEE\,$ the total manifold
and $\,\MMM\,$ the base manifold,
i.e. $\,\Range{\Bp}=\MMM\,$ and
$\,\Range{\TT\Bp\di\Be}=\TANG_{\Bp\di\Be}\MMM\,$ for all $\,\Be\in\EEE\,$.
The vertical distribution is $\,\VE\equaldef\Null{\TT\Bp}\,$.
A section
$\,\Bs\in\cont^1\di{\MMM\sp\EEE}\,$ is such that
$\,\Bp\circ\Bs\in\cont^1\di{\MMM\sp\MMM}\,$ is the identity.
The fibre at $\,\Bx\in\MMM\,$ is the set 
$\,\EEE_\Bx\equaldef\Bp^{-1}\di\Bx\,$
which is assumed to be isomorphic to a standard fibre manifold.
The pull-back bundle of the tangent bundle $\,\TprojE\in\cont^1\di{\TE\sp\EEE}\,$
by a section
$\,\Bs\in\cont^1\di{\MMM\sp\EEE}\,$ is the fibre bundle
$\,\Bs\pull\TprojE\in\cont^1\di{\Bs\pull\TE\sp\MMM}\,$
whose fibre at $\,\Bx\in\MMM\,$
is the tangent space $\,\TANG_{\Bs\di\Bx}\EEE\,$ of 
$\,\TprojE\in\cont^1\di{\TE\sp\EEE}\,$.
\begin{definition}[Connection]\label{def: connection}
\index{connection}\index{horizontal lift}
A connection
$\,\VPROJ\in\Lambda^1\di{\EEE\sp\TE}\,$
in a fibre bundle
$\,\Bp\in\cont^1\di{\EEE\sp\MMM}\,$
is an idempotent vector-valued one-form,
which is pointwise a projector on vertical subspaces:
$\,\VPROJ\circ\VPROJ=\VPROJ\,$ with
$\,\Range{\VPROJ\di\Be}=\Null{\TT\Bp\di\Be}\,$.
Horizontal vectors are the ones
in the kernel $\,\Null{\VPROJ\di\Be}\,$ of the connection.
The projector on the
horizontal distribution $\,\HE\,$ is denoted by
$\,\HPROJ=\ID\TE-\VPROJ\,$, so that
$\,\HPROJ\circ\HPROJ=\HPROJ\,$ and
$\,\HPROJ\circ\VPROJ=\VPROJ\circ\HPROJ=0\,$.
\end{definition}
The tangent to a section
$\,\Bs\in\cont^1\di{\MMM\sp\EEE}\,$
of a fibre bundle
$\,\Bp\in\cont^1\di{\EEE\sp\MMM}\,$
along a vector field
$\,\Bv\in\cont^0\di{\MMM\sp\TM}\,$
is a section
$\,\TT\Bs\punto\Bv\in\cont^1\di{\MMM\sp\Bs\pull\TE}\,$
of the pull-back bundle $\,\Bs\pull\Bp=\cont^1\di{\Bs\pull\TE\sp\MMM}\,$.
\begin{definition}[Natural derivative]
\label{def: naturalder}\index{natural derivative}
In a fibre bundle
$\,\Bp\in\cont^1\di{\EEE\sp\MMM}\,$,
the natural derivative of a section
$\,\Bs\in\cont^1\di{\MMM\sp\EEE}\,$
according to a vector field
$\,\Bv\in\cont^0\di{\MMM\sp\TM}\,$
is the tangent vector field
$\,\Tder\Bv\in\cont^1\di{\Bs\di\MMM\sp\TE}\,$
in the tangent bundle 
$\,\TprojE\in\cont^1\di{\TE\sp\EEE}\,$
defined by
$$\Tder\Bv\circ\Bs\equaldef\TT\Bs\punto\Bv\in\cont^1\di{\MMM\sp\TE}\,.$$
\end{definition}
For any $\,\Bx\in\MMM\,$ we have that
$\,\Tder\Bv\di{\Bs_\Bx}=\TT_{\Bv_\Bx}\Bs\in\TANG_{\Bs_\Bx}\EEE\,$.
The natural derivative 
$\,\Tder\Bv\in\cont^1\di{\Bs\di\MMM\sp\TE}\,$
is $\,\Bp$-related to the vector field
$\,\Bv\in\cont^1\di{\MMM\sp\TM}\,$, because:
$$\vcenter{\halign{
\hfil$#$&$#$\hfil&$#$\hfil&$#$\hfil\cr
\TT\Bp\circ\Tder\Bv&\,=\Bv\circ\Bp\in\cont^1\di{\Bs\di\MMM\sp\TM}\,,
\vspace{8pt}\cr
\Bp\circ\FL{\Tder\Bv}\pttt&\,=
\FL\Bv\pttt\circ\Bp\in\cont^1\di{\EEE\sp\MMM}
\,.\cr}}$$
It is also apparent that the natural derivative is tensorial in
$\,\Bv\in\cont^0\di{\MMM\sp\TM}\,$ since
the differential $\,\Tder{\Bv_\Bx}\Bs\in\TANG_{\Bs\di\Bx}\EEE\,$
is linearly dependent on the vector
$\,\Bv_\Bx\in\TANG_{\Bx}\MMM\,$.
The next statement enunciates a well known property of naturality of the \Lie\ bracket
with respect to relatedness, 
(see e.g. 
\cite{KMS1993} Lemma 3.10
or \cite{RomanoDiff2007} Lemma 1.3.4).
\begin{lemma}[Morphism-related vector fields and Lie brackets]
\label{lm: relatedbrackets}
Let the vector fields 
$\,\BX,\BY\in\cont^1\di{\NNN\sp\TN}\,$ 
be related to vector fields
$\,\Bu,\Bv\in\cont^1\di{\MMM\sp\TM}\,$ 
by a morphism 
$\,\map\in\cont^1\di{\MMM\sp\NNN}\,$, viz:
$$\vcenter{\halign{
\hfil$#$&$#$\hfil&$#$\hfil&$#$\hfil\cr
\BX\circ\map
&\,=\TT\map\circ\Bu\,,
\qquad
\BY\circ\map
&\,=\TT\map\circ\Bv
\,.\cr}}$$
Then also their \Lie\ brackets are $\,\map$-related:
$$\,\lineare{\BX,\BY}\circ\map
=\TT\map\circ\lineare{\Bu,\Bv}\,.$$
Setting $\,\Tder\Bv\circ\map\equaldef\TT\map\circ\Bv\,$ 
for any morphism 
$\,\map\in\cont^1\di{\MMM\sp\NNN}\,$, we have that
$\,\TT\map\circ\lineare{\Bu,\Bv}=\Tder{\lineare{\Bu,\Bv}}\circ\map\,$
and the result may be stated as
$\,\lineare{\Tder\Bu,\Tder\Bv}
=\Tder{\lineare{\Bu,\Bv}}\,$.
\end{lemma}
Tensoriality is a crucial property of a multilinear 
scalar or vector valued map, 
meaning that it \textit{lives at points} \cite{Spivak1979}, i.e. that
its point-values depend only on the values of the argument fields at that point.
A standard tensoriality criterion for multilinear forms
on $\,\MMM\,$ is provided by $\,\cont^\infty\di{\MMM\sp\Re}$-linearity 
(see \cite{KMS1993} Lemma 7.3 
or \cite{Lang1995} Lemma 2.3 of Ch. VIII).

Although not needed in evaluating the \Lie\ bracket
$\brack{\Tder\Bu}{\Tder\Bv}$ on $\,\Bs\di\MMM\,$,
for the developments illustrated in Theorem \ref{th: curvnice}
it is essential to extend the domain of the natural derivatives 
$\,\Tder\Bu,\Tder\Bv\in\cont^1\di{\Bs\di\MMM\sp\TE}\,$
outside the range $\,\Bs\di\MMM\subset\EEE\,$
of the section $\,\Bs\in\cont^1\di{\MMM\sp\EEE}\,$,
so that they can be considered as
(local) tangent vector fields 
$\,\Tder\Bu,\Tder\Bv\in\cont^1\di{\EEE\sp\TE}\,$
with the further property of being projectable.
This task can be accomplished by the following construction.
\begin{lemma}[Extension by foliation]\label{lm: extfoliation}
The natural derivative of a section
$\,\Bs\in\cont^1\di{\MMM\sp\EEE}\,$
of a fibre bundle
$\,\Bp\in\cont^1\di{\EEE\sp\MMM}\,$,
according to a vector field
$\,\Bv\in\cont^0\di{\MMM\sp\TM}\,$,
can be extended,
in the bundle 
$\,\TprojE\in\cont^1\di{\TE\sp\EEE}\,$,
to a (local) tangent vector field
$\,\Tder\Bv\in\cont^1\di{\EEE\sp\TE}\,$
which projects on the vector field
$\,\Bv\in\cont^0\di{\MMM\sp\TM}\,$, i.e.
we have that, locally in $\,\EEE\,$:
$$\vcenter{\halign{
$#$\hfil&$#$\hfil&$#$\hfil&$#$\hfil&$#$\hfil\cr
&\,\TprojE\circ\Tder\Bv=\ID{\EEE}\,,
\vspace{8pt}\cr
&\,\TT\Bp\circ\Tder\Bv=\Bv\circ\Bp
\,.\cr}}$$
\end{lemma}
\proof
The extension may be performed by considering
a (local) foliation of the total manifold $\,\EEE\,$,
whose leaves are transversal to the fibres and
include the folium $\,\Bs\di\MMM\,$.
The existence of at least a local foliation with these characteristics
can be inferred by acting with a local bundle chart, 
which maps (locally) the image of the section into the 
trivial bundle image of the chart, 
and, subsequently, with a local chart which maps 
(locally) the fibres in their linear model space.
The foliation is then performed by translation
in the linear image of the fibres and the resulting leaves are mapped 
back to get the leaves in the total manifold.
It is thus possible to define the map
$\,\Bss\in\cont^1\di{\EEE\sp\cont^1\di{\MMM\sp\EEE}}\,$
which to each $\,\Be\in\EEE\,$ associates the (local)
section $\,\Bss_\Be\in\cont^1\di{\MMM\sp\EEE}\,$ by
$$\Bss_\Be\di\Bx\equaldef
\BSigma_\Be\cap\EEE_\Bx
\,,\perogni\Be\in\EEE\,,$$
whose range is the leaf $\,\BSigma_\Be\,$ through $\,\Be\in\EEE\,$.
The extension of $\,\Tder\Bv\,$
is then (locally) defined by
$\,\Tder\Bv\di\Be\equaldef
\TT_{\Bp\di\Be}\Bss_\Be\punto\Bv_{\Bp\di\Be}\,$
and gives a vector field since
$\,\TprojE\di{\TT_{\Bp\di\Be}\Bss_\Be\punto\Bv_{\Bp\di\Be}}=\Be\,$
for all $\,\Be\in\EEE\,$.
Moreover this extension projects on $\,\Bv\in\cont^0\di{\MMM\sp\TM}\,$ since
$$\vcenter{\halign{
$#$\hfil&$#$\hfil&$#$\hfil\cr
\TT_{\Bp\di\Be}\Bp\punto\Tder\Bv\di\Be
&\,=\TT_{\Bp\di\Be}\Bp\punto\TT_{\Bp\di\Be}\Bss_\Be\punto\Bv_{\Bp\di\Be}
=\TT_{\Bp\di\Be}(\Bp\circ\Bss_\Be)\punto\Bv_{\Bp\di\Be}
=\Bv_{\Bp\di\Be}
\,.\cr}}$$
Being
$\,\Bss_\Be\di{\Bp\di\Be}=\Be\,$
the extension
$\,\Tder\Bv\di\Be\equaldef
\TT_{\Bp\di\Be}\Bss_\Be\punto\Bv_{\Bp\di\Be}\,$
may be written as
$\,(\Tder\Bv\circ\Bss_\Be)\di{\Bp\di\Be}
=(\TT\Bss_\Be\circ\Bv)\di{\Bp\di\Be}\,$
which, by surjectivity of $\,\Bp\,$, means that (locally)
$$\Tder\Bv\circ\Bss_\Be=\TT\Bss_\Be\circ\Bv
\,,\perogni\Bx\in\MMM\,.$$
If $\,\Be_1,\Be_2\in\EEE\,$ are such that $\,\BSigma_{\Be_1}=\BSigma_{\Be_2}\,$,
then $\,\Bss_{\Be_1}=\Bss_{\Be_2}\,$.
If $\,\Be\in\Bs\di\MMM\,$, the section 
$\,\Bss_\Be\in\cont^1\di{\MMM\sp\EEE}\,$
is in fact coincident with
$\,\Bs\in\cont^1\di{\MMM\sp\EEE}\,$.
\endprova
\begin{definition}[Horizontal lift]\index{horizontal lift}
In a bundle
$\Bp\in\cont^1\di{\EEE\sp\MMM}$
the horizontal lift 
$\,\Hlift\in\cont^1\di{\EEE\times_\MMM\TM\sp\TE}\,$
is a right inverse of 
$\,\coppia{\TprojE}{\TT\Bp}\in\cont^1\di{\TE\sp\EEE\times_\MMM\TM}\,$
such that the map
$\,\Hlift_{\Bs_\Bx}\in\cont^1\di{\TM\sp\TE}\,$,
defined by
$\,\Hlift_{\Bs_\Bx}\di{\Bv_\Bx}=\Hlift\coppia{\Bs_\Bx}{\Bv_\Bx}\,$
for all $\,\Bv_\Bx\in\TANG_\Bx\MMM\,$,
is a linear homomorphism from the tangent bundle
$\,\TprojM\in\cont^1\di{\TM\sp\MMM}\,$ to the tangent bundle
$\,\TprojE\in\cont^1\di{\TE\sp\EEE}\,$, i.e:
$$\vcenter{\halign{
\hfil$#$&$#$\hfil&$#$\hfil\cr
&\,\coppia{\TprojE}{\TT\Bp}\circ\Hlift=\ID{\EEE\times_\MMM\TM}\,,
\vspace{6pt}\cr
&\,\Hlift_{\Bs_\Bx}\di{\aa\,\Bu_\Bx+\bb\,\Bv_\Bx}
=\aa\,\Hlift_{\Bs_\Bx}\di{\Bu_\Bx}+\bb\,\Hlift_{\Bs_\Bx}\di{\Bv_\Bx}
\in\TANG_{\Bs_\Bx}\EEE\,,\cr}}$$
with $\,\Bs_\Bx\in\EEE_\Bx\,$ and $\,\Bu_\Bx,\Bv_\Bx\in\TANG_\Bx\MMM\,$
and $\,\aa,\bb\in\Re\,$.
\end{definition}
\begin{lemma}[Horizontal lifts and horizontal projectors]
\label{lm: horlifthorproj}
Given a horizontal projector $\,\HPROJ\in\cont^1\di{\TE\sp\TE}\,$,
the induced horizontal lift is defined by
$$\Hlift\coppia{\Bs_\Bx}{\Bv_\Bx}\equaldef
\HPROJ\punto\TT_\Bx\Bs\punto\Bv_\Bx\in\HORI_{\Bs_\Bx}\BUNDLE
\,,\perogni\Bs_\Bx\in\EEE_\Bx,\quad\Bv_\Bx\in\TANG_\Bx\MMM\,,$$
where $\,\Bs\in\cont^1\di{\MMM\sp\EEE}\,$ is an arbitrary section
extension of $\,\Bs_\Bx\in\EEE_\Bx\,$.
Vice versa, a horizontal lift
$\,\Hlift\in\cont^1\di{\EEE\times_\MMM\TM\sp\TE}\,$
induces a horizontal projector given by
$\,\HPROJ\equaldef\Hlift\circ\coppia{\TprojE}{\TT\Bp}\,$.
\end{lemma}
\proof
The former formula yields a horizontal lift since:
$$(\coppia{\TprojE}{\TT\Bp}\circ\Hlift)\coppia{\Bs_\Bx}{\Bv_\Bx}
=\coppia{\TprojE}{\TT\Bp}\punto\HPROJ\punto\TT_\Bx\Bs\punto\Bv_\Bx
=\coppia{\Bs_\Bx}{\Bv_\Bx}\,,$$
and the latter formula yields a horizontal projector because the homomorphism
$\,\HPROJ\equaldef\Hlift\circ\coppia{\TprojE}{\TT\Bp}\,$
is idempotent by
$\,\HPROJ\circ\HPROJ=
\Hlift\circ\coppia{\TprojE}{\TT\Bp}
\circ
\Hlift\circ\coppia{\TprojE}{\TT\Bp}
=\Hlift\circ\ID{\EEE\times_\MMM\TM}\circ\coppia{\TprojE}{\TT\Bp}=\HPROJ
\,$
and horizontal by the identity
$\,(\coppia{\TprojE}{\TT\Bp}\circ\HPROJ)\di\BX
=(\coppia{\TprojE}{\TT\Bp}\circ\Hlift\circ\coppia{\TprojE}{\TT\Bp})\di\BX
=\coppia{\TprojE\di\BX}{\TT\Bp\di\BX}\,$.
\endprova
\begin{definition}[Covariant derivative]
\label{def: Horicovar}
The covariant derivative is the vertical component of the natural derivative:
$$\vcenter{\halign{
\hfil$#$&$#$\hfil&$#$\hfil\cr
\nablabar_\Bv\Bs&\,\equaldef\VPROJ\circ\Tder\Bv\circ\Bs
\in\cont^1\di{\MMM\sp\VE}
\,.\cr}}$$
\end{definition}
Setting $\,\Hlift\Bs=\HPROJ\circ\TT\Bs\,$
and $\,\nablabar\Bs=\VPROJ\circ\TT\Bs\,$,
it is
$\,\TT\Bs=\nablabar\Bs+\Hlift\Bs
\in\cont^1\di{\TM\sp\TE}\,$
and
$\,\TT_\Bv=\nablabar_\Bv+\Hlift_\Bv\in\cont^1\di{\Bs\di\MMM\sp\TE}\,$
with $\,\nablabar_\Bv=\VPROJ\circ\TT_\Bv\,$ and $\,\Hlift_\Bv=\HPROJ\circ\TT_\Bv\,$.
\begin{lemma}[Projectability]
\label{lm: horirelat}
The horizontal lift
$\,\Hlift_\Bv\in\cont^1\di{\Bs\di\MMM\sp\HE}\,$
is $\,\Bp$-related to the vector field
$\,\Bv\in\cont^{1}\di{\MMM\sp\TM}\,$:
$\,\TT\Bp\circ\Hlift_\Bv=\Bv\circ\Bp\in\cont^0\di{\Bs\di\MMM\sp\TM}\,$.
\end{lemma}
\proof
From the decomposition
$\,\TT_\Bv=\nablabar\Bv+\Hlift_\Bv\in\cont^1\di{\EEE\sp\TE}\,$
it follows that:
$\,\TT\Bp\circ\TT_\Bv
=\TT\Bp\circ\nablabar_\Bv+\TT\Bp\circ\Hlift_\Bv
=\TT\Bp\circ\Hlift_\Bv\,$
being, by definition $\,\TT\Bp\circ \nablabar_\Bv=0\,$.
The $\,\Bp$-relatedness of $\,\Hlift_\Bv\,$ to $\,\Bv\,$
is then inferred from that of $\,\TT_\Bv\,$.
\endprova
Naturality of \Lie\ brackets with respect to relatedness and
Lemma \ref{lm: horirelat} give:
$$\TT\Bp\circ\brack{\Hlift_\Bu}{\Hlift_\Bv}
=\brack{\TT\Bp\circ\Hlift_\Bu}{\TT\Bp\circ\Hlift_\Bv}
=\brack{\Bu\circ\Bp}{\Bv\circ\Bp}
=\brack{\Bu}{\Bv}\circ\Bp
\in\cont^1\di{\EEE\sp\TM}\,.$$%
\begin{lemma}[Injectivity]
\label{lm: horizontallift}
The horizontal lift
$\,\Hlift\Bs\in\cont^1\di{\TM\sp\HE}\,$,
along a cross section $\,\Bs\in\cont^1\di{\MMM\sp\EEE}\,$
of a fibre bundle
$\,\Bp\in\cont^1\di{\EEE\sp\MMM}\,$,
is a fibrewise injective homomorphism, i.e.
$\,\Hlift_\Bx\Bs\in\Linmap{\TANG_\Bx\MMM\sp\HORI_{\Bs\di\Bx}\EEE}\,$
is an injective linear map at each $\,\Bx\in\MMM\,$.
\end{lemma}
\proof
We must prove that $\,\Null{\Hlift_\Bx\Bs}=\set0\,$.
We first investigate the linear differential
$\,\TT_\Bx\Bs\in\Linmap{\TANG_\Bx\MMM\sp\TANG_{\Bs\di\Bx}\EEE}\,$.
By the characteristic property of a section,
$\,\Bproj\circ\Bs=\ID\MMM\,$ it is:
$\,\TT_{\Bs\di\Bx}\Bproj\punto\TT_\Bx\Bs\punto\Bv_\Bx
=\TT_\Bx(\Bproj\circ\Bs)\punto\Bv_\Bx
=\Bv_\Bx\,$ for all
$\,\Bv_\Bx\in\TANG_\Bx\MMM\,$.
It follows that $\,\Null{\TT_\Bx\Bs}=\set0\,$
and $\,\Range{\TT_\Bx\Bs}\cap\Null{\TT_{\Bs\di\Bx}\Bproj}=\set0\,$.
The injectivity of $\,\TT_\Bx\Bs\,$ implies that:
$\,\dim\Range{\TT_\Bx\Bs}=\dim\TANG_\Bx\MMM\,$.
Being
$\,\TT_\Bx\Bs=\nablabar_\Bx\Bs+\Hlift_\Bx\Bs\,$ with
$\,\Range{\nablabar_\Bx\Bs}\subseteq\Null{\TT_{\Bs\di\Bx}\Bproj}\,$,
we have that
$\,\TT_{\Bs\di\Bx}\Bproj\punto\Hlift_\Bx\Bs\punto\Bv_\Bx
=\TT_{\Bs\di\Bx}\Bproj\punto\TT_\Bx\Bs\punto\Bv_\Bx
=\Bv_\Bx\,$ for all
$\,\Bv_\Bx\in\TANG_\Bx\MMM\,$. 
It follows that  $\,\Null{\Hlift_\Bx\Bs}=\set0\,$
and
$\,\Range{\Hlift_\Bx\Bs}\cap\Null{\TT_{\Bs\di\Bx}\Bproj}=\set0\,$
with
$\,\dim\Range{\Hlift_\Bx\Bs}=\dim\TANG_\Bx\MMM\,$.
\endprova
\begin{theorem}[Homomorphism]
\label{th: horizontallift}
The horizontal lift
$\,\Hlift\Bs\in\cont^1\di{\TM\sp\HE}\,$
along a section $\,\Bs\in\cont^1\di{\MMM\sp\EEE}\,$
of a fibre bundle
$\,\Bp\in\cont^1\di{\EEE\sp\MMM}\,$
is a vector bundle homomorphism between the bundle
$\,\TprojM\in\cont^1\di{\TM\sp\MMM}\,$ and the pull-back bundle 
$\,\Bs\pull\TprojE\in\cont^1\di{\Bs\pull\HE\sp\MMM}\,$
which is fibrewise invertible and
tensorial in $\,\Bs\in\cont^1\di{\MMM\sp\EEE}\,$.
\end{theorem}
\proof
Let $\,\dim\MMM=\dim\TANG_\Bx\MMM=m\,$ and $\,\dim\FFF=f\,$ 
where $\,\FFF\,$ is the typical fibre.
Then
$\,\dim\EEE=\dim\TANG_{\Bs\di\Bx}\EEE=m+f\,$. 
So that $\,\dim\VERT_{\Bs\di\Bx}\EEE=f\,$ and
$\,\dim\HORI_{\Bs\di\Bx}\EEE=m\,$. 
By reasons of dimensions the injectivity of
$\,\Hlift_\Bx\Bs\in\Linmap{\TANG_\Bx\MMM\sp\HORI_{\Bs\di\Bx}\EEE}\,$ 
implies then its surjectivity.
Moreover let
$\,\overline\Bs\in\cont^1\di{\MMM\sp\EEE}\,$
be another section such that
$\,\overline\Bs\di\Bx=\Bs\di\Bx\,$.
Then, for any $\,\Bv_\Bx\in\TANG_\Bx\MMM\,$, being 
$\,\TT_{\Bv_\Bx}\Bs,\TT_{\Bv_\Bx}\overline\Bs\in\TANG_{\Bs\di\Bx}\EEE\,$,
we have that
$\,\TT\Bp\circ(\TT_{\Bv_\Bx}\Bs-\TT_{\Bv_\Bx}\overline\Bs)=0\,$
and hence that 
$\,\Hlift_{\Bv_\Bx}\Bs=\HPROJ\circ\TT_{\Bv_\Bx}\Bs
=\HPROJ\circ\TT_{\Bv_\Bx}\overline\Bs=\Hlift_{\Bv_\Bx}\overline\Bs
\in\Linmap{\TANG_\Bx\MMM\sp\HORI_{\Bs\di\Bx}\EEE}\,$.
To a tangent vector $\,\Bv_\Bx\in\TANG_\Bx\MMM\,$ 
there corresponds a horizontal vector
$\,\Hlift_{\Bv_\Bx}\Bs\in\HORI_{\Bs\di\Bx}\EEE\,$
which depends only on the value of
$\,\Bs\in\cont^1\di{\MMM\sp\EEE}\,$ at $\,\Bx\in\MMM\,$.
\endprova

\section{Curvature of a connection}
\label{sec: Curvature}

The vertical distribution of a fibre bundle
$\,\Bp\in\cont^1\di{\EEE\sp\MMM}\,$ is integrable and
the leaves of the induced foliation are the fibres of the bundle.
By \Frobenius\ theorem \cite{KMS1993}, \cite{Lang1995},
integrability of vertical distribution is inferred form 
the vanishing of the vector-valued \textit{cocurvature} form:
$\,\cocurvform\di{\BX,\BY}
\equaldef-\HPROJ\circ\lineare{\widehat{\VPROJ\BX},\widehat{\VPROJ\BY}}=0\,$
for any $\,\BX,\BY\in\TE\,$.
Here
$\,\coppia{\widehat{\VPROJ\BX}}{\widehat{\VPROJ\BY}}\in\cont^1\di{\EEE\sp\TE}\,$
is any pair of vector fields extension of the vectors 
$\, \VPROJ\BX,\VPROJ\BY\in\TE\,$, since  tensoriality
follows from the 
$\,\cont^\infty\di{\EEE\sp\Re}$-linearity of the cocurvature form.
The involutivity condition:
$\brack{\widehat{\HPROJ\BX}}{\widehat{\HPROJ\BY}}\in\cont^1\di{\EEE\sp\HE}\,$,
to be imposed for the integrability of the horizontal distribution,
is equivalently expressed by 
the vanishing of the \textit{curvature}
defined by \cite{KMS1993}:\index{curvature}
$$\curvform\di{\BX,\BY}\equaldef
-\VPROJ\circ\brack{\widehat{\HPROJ\BX}}{\widehat{\HPROJ\BY}}\,,
\perogni\BX,\BY\in\TE\,.$$
Again tensoriality follows from the 
$\,\cont^\infty\di{\EEE\sp\Re}$-linearity of the curvature form,
as shown below.
Let us denote by $\,\Lambda^k\di{\MMM\sp\TM}\,$
the space of tangent-valued $\,k$-forms on a manifold $\,\MMM\,$.
\begin{proposition}[Tensoriality of the curvature]
\label{prop: curvtens}
The curvature of a connection 
$\,\VPROJ\in\Lambda^1\di{\EEE\sp\TE}\,$
in a fibre bundle
$\,\Bproj\in\cont^1\di{\EEE\sp\MMM}\,$
is a vertical-vector valued, horizontal $\,2$-form
$\,\curvform\in\Lambda^2\di{\EEE\sp\VE}\,$, that is
a $\,2$-form vanishing on
vertical vectors and taking values in the vertical distribution. 
\end{proposition}
\proof
A direct verification of the tensoriality, based on 
$\,\cont^\infty\di{\EEE\sp\Re}$-linearity, yields the result:
$$\vcenter{\halign{
\hfil$#$&$#$\hfil&$#$\hfil\cr
-\curvform\di{\BX,\ff\BY}\,\equaldef&\,
\VPROJ\circ\brack{\widehat{\HPROJ\BX}}{\widehat{\HPROJ\BY}}
\vspace{8pt}\cr
=&\,\ff\,\VPROJ\circ\brack{\widehat{\HPROJ\BX}}{\widehat{\HPROJ\BY}}
+(\Lieder_{\HPROJ\BX}\ff)\,(\VPROJ\circ\HPROJ)\di\BY
\vspace{8pt}\cr
=&\,-\ff\,\curvform\di{\BX,\BY}\,,
\perogni\ff\in\cont^1\di{\EEE\sp\Re}
\,,\cr}}$$
since $\,\VPROJ\circ\HPROJ=0\,$.
Similarly $\,\curvform\di{\ff\BX,\BY}=\ff\,\curvform\di{\BX,\BY}\,$.
\endprova
\begin{theorem}
\label{th: curvhori}
For any given section 
$\,\Bs\in\cont^1\di{\MMM\sp\EEE}\,$,
the curvature of a connection
$\,\VPROJ\in\Lambda^1\di{\EEE\sp\VE}\,$
is expressed by a $\,2$-form
$\,\curvoperator_\Bs\in\Lambda^2\di{\MMM\sp\Bs\pull\VE}\,$
with values in the pull-back of the vertical distribution 
by the section 
$\,\Bs\in\cont^1\di{\MMM\sp\EEE}\,$,
defined in terms of horizontal lifts by:
$$\curvopbar_\Bs\di{\Bu,\Bv}\equaldef
\curvform\di{\Hlift_\Bu,\Hlift_\Bv}\circ\Bs
=(\Hlift_{\brack\Bu\Bv}-\brack{\Hlift_\Bu}{\Hlift_\Bv})\circ\Bs\,,
\perogni\Bu,\Bv\in\Lambda^0\di{\MMM\sp\TM}\,,$$
The $\,2$-form 
$\,\curvoperator_\Bs\in\Lambda^2\di{\MMM\sp\Bs\pull\VE}\,$
is tensorial in $\,\Bs\in\cont^1\di{\MMM\sp\EEE}\,$.
\end{theorem}
\proof
We rely on the properties of tensoriality and horizontality
of the curvature two-form
$\,\curvform\in\Lambda^2\di{\EEE\sp\VE}\,$
stated in Proposition \ref{prop: curvtens}
and on the tensorial isomorphism of the horizontal lifts
stated in Theorem \ref{th: horizontallift}.
Accordingly, the point value of the curvature 
$\,\curvform\di{\BX,\BY}=
-\VPROJ\circ\brack{\widehat{\HPROJ\BX}}{\widehat{\HPROJ\BY}}\,$ 
at $\,\Bb\in\EEE_\Bx\,$
depends only on the vectors 
$\,\HPROJ\BX_\Bb,\HPROJ\BY_\Bb\in\TANG_\Bb\EEE\,$.
Moreover, 
by Theorem \ref{th: horizontallift},
given any section $\,\Bs\in\cont^1\di{\MMM\sp\EEE}\,$
such that $\,\Bs_\Bx=\Bb\,$, 
there exists a uniquely determined pair of vectors
$\,\Bu_\Bx,\Bv_\Bx\in\TANG_\Bx\MMM\,$, such that
$\,\Hlift_{\Bu_\Bx}\Bs=(\HPROJ\BX)\di{\Bs_\Bx}\,,\quad
\Hlift_{\Bv_\Bx}\Bs=(\HPROJ\BY)\di{\Bs_\Bx}\,$
and the pair $\,\Bu_\Bx,\Bv_\Bx\in\TANG_\Bx\MMM\,$ does not
depend on the choice of the section $\,\Bs\in\cont^1\di{\MMM\sp\EEE}\,$
such that $\,\Bs_\Bx=\Bb\,$.
Then the curvature two-form
$\,\curvform\in\Lambda^2\di{\EEE\sp\VE}\,$,
evaluated on pairs of horizontal lifts,
defines the field
$\,\curvoperator_{\Bs}\di{\Bu,\Bv}
\equaldef
-\VPROJ\circ\lineare{\Hlift_\Bu,\Hlift_\Bv}\circ\Bs
\in\cont^1\di{\MMM\sp\VE}\,$
for any pair of vector fields
$\,\Bu,\Bv\in\cont^0\di{\MMM\sp\TM}\,$
on the tangent bundle and any section
$\,\Bs\in\cont^1\di{\MMM\sp\EEE}\,$
of the fibre bundle
$\,\Bproj\in\cont^1\di{\EEE\sp\MMM}\,$.
By tensoriality, for any section
$\,\Bs\in\cont^1\di{\MMM\sp\EEE}\,$ the field 
$\,\curvoperator_\Bs\in\Lambda^2\di{\MMM\sp\VE}\,$
is a vector-valued two-form
on $\,\MMM\,$ with values in $\,\Bs\pull\VE\,$
and for any pair
$\,\Bu,\Bv\in\cont^0\di{\MMM\sp\TM}\,$
the field 
$\,\curvoperator\di{\Bu,\Bv}\in\Lambda^1\di{\MMM\sp\Bs\pull\VE}\,$
is a vertical-valued vector field
along $\,\Bs\in\cont^1\di{\MMM\sp\EEE}\,$.
Moreover, by Lemma \ref{lm: horirelat}, the horizontal lifts are 
projectable and we have the relations:
$$\left.\vcenter{\halign{
\hfil$#$&$#$\hfil&$#$\hfil\cr
\TT\Bp\circ\brack{\Hlift_\Bu}{\Hlift_\Bv}
&\,=\brack{\Bu}{\Bv}\circ\Bp
\vspace{6pt}\cr
\TT\Bp\circ\Hlift_{\brack{\Bu}{\Bv}}
&\,=\brack{\Bu}{\Bv}\circ\Bp
\cr}}\right\}
\implies
\TT\Bp\circ(\lineare{\Hlift_\Bu,\Hlift_\Bv}-\Hlift_{\lineare{\Bu,\Bv}})=0
\,.$$
Then $\,\Hlift_{\brack\Bu\Bv}\,$ 
is the horizontal component of 
$\,\brack{\Hlift_\Bu}{\Hlift_\Bv}\,$ and we get the equality:
$\,\brack{\Hlift_\Bu}{\Hlift_\Bv}-\Hlift_{\brack\Bu\Bv}=
\VPROJ\circ\brack{\Hlift_\Bu}{\Hlift_\Bv}
\equi\Hlift_{\brack\Bu\Bv}=\HPROJ\circ\brack{\Hlift_\Bu}{\Hlift_\Bv}\,$.
\endprova

\section{Covariant derivative}
\label{sec: covariant}

\begin{lemma}[Covariant derivative as Lie derivative]
\label{lm: covarLie}
In a fibre bundle
$\,\Bproj\in\cont^1\di{\EEE\sp\MMM}\,$ 
with a connection, the covariant deri\-vative 
may be defined as the generalized \Lie\ derivative:
$$\vcenter{\halign{
\hfil$#$&$#$\hfil&$#$\hfil&$#$\hfil&$#$\hfil\cr
\nablabar_\Bv\Bs
=\Lieder_{\coppia{\Hlift_\Bv}{\Bv}}\Bs
=\parder{\pttt}{0}\FL{\coppia{\Hlift_\Bv}{\Bv}}{\pttt}\pull\Bs
=\parder{\pttt}{0}
\FL{\Hlift_\Bv}{-\pttt}\circ\Bs\circ\FL{\Bv}{\pttt}\,.
\cr}}$$
\end{lemma}
\proof
By \Leibniz\ rule
$\Lieder_{\coppia{\Hlift_\Bv}{\Bv}}\Bs
=\TT\Bs\circ\Bv-\Hlift_\Bv\Bs
=\TT_\Bv\Bs-\Hlift_\Bv\Bs\,$.
Then, being 
$\,\Lieder_{\coppia{\Hlift_\Bv}{\Bv}}\Bs\in\cont^1\di{\MMM\sp\VE}\,$
and 
$\,\Hlift_\Bv\Bs\in\cont^1\di{\MMM\sp\HE}\,$,
by uniqueness of the vertical-horizontal split, we get that 
$\,\nablabar_\Bv\Bs\equaldef\VPROJ\circ\Tder\Bv\Bs
=\Lieder_{\coppia{\Hlift_\Bv}{\Bv}}\Bs\,$.
\endprova
\begin{definition}[Parallel transport]
\label{def: parallel}\index{parallel transport}
Let $\,\Bproj\in\cont^1\di{\EEE\sp\MMM}\,$
be a fibre bundle with a connection.
The parallel transport 
$\,\FL{\Bv}{\pttt}\forw\Bs\in\cont^1\di{\MMM\sp\EEE}\,$
of a section
$\,\Bs\in\cont^1\di{\MMM\sp\EEE}\,$
along the flow 
$\,\FL{\Bv}{\pttt}\in\cont^1\di{\MMM\sp\MMM}\,$
is defined by:
$$\FL{\Bv}{\pttt}\forw\Bs\equaldef
\FL{\Hlift_\Bv}{\pttt}\circ\Bs
=(\FL{\coppia{\Hlift_\Bv}{\Bv}}{\pttt}\push\Bs)\circ\FL{\Bv}{\pttt}\,,$$ 
so that
$\,\Bp\circ\FL{\Bv}{\pttt}\forw\Bs
=\Bp\circ\FL{\Hlift_\Bv}{\pttt}\circ\Bs
=\FL{\Bv}{\pttt}\circ\Bp\circ\Bs
=\FL{\Bv}{\pttt}\,$.
\end{definition} 
From the definition of parallel transport
and Lemma \ref{lm: covarLie} we infer
that the covariant derivative 
and the horizontal lift are given by:
$$\vcenter{\halign{
\hfil$#$&$#$\hfil&$#$\hfil\cr
&\,\nablabar_\Bv\Bs
=\parder{\pttt}{0}\FL{\Hlift_\Bv}{-\pttt}\circ\Bs\circ\FL{\Bv}{\pttt}
=\parder{\pttt}{0}\FL{\Bv}{-\pttt}\forw\Bs\circ\FL{\Bv}{\pttt}\,,
\vspace{8pt}\cr
&\,\Hlift_\Bv\Bs
=\parder{\pttt}{0}\FL{\Hlift_\Bv}{\pttt}\circ\Bs
=\parder{\pttt}{0}\FL{\Bv}{\pttt}\forw\Bs
\,.\cr}}$$
Since the horizontal lift $\,\Hlift_\Bv\,$
is defined pointwise in $\,\MMM\,$,
the parallel transport along a curve in $\,\MMM\,$
of a section defined only on that curve is meaningful
and so is for the covariant derivative.
\begin{definition}[Geodesic]\index{geodesic curve}
A curve $\,\cc\in\cont^1\di{\II\sp\MMM}\,$ in a manifold with
a connection is a geodesic if 
the velocity field of the curve
$\,\Bv\in\cont^1\di{\II\sp\TM}\,$
fulfils the condition
$$\nablabar_{t}\Bv\equaldef\parder\tau\ttt\cc_{\ttt,\tau}\forw\Bvtau=0\,,$$
where $\,\nablabar\,$ is the covariant derivative,
the velocity is given by $\,\Bv_\ttt\equaldef\parder\tau\ttt\cc_\tau\,$
and $\,\cc_{\ttt,\tau}\forw\,$ is the parallel transport from $\,\cc_\tau\,$
to $\,\cc_\ttt\,$ along the curve.
\end{definition}
\begin{definition}[Spray]\label{def: flipdue}
A section $\,\BX\in\cont^1\di{\TM\sp\TTM}\,$
of the tangent bundle $\,\Tproj_\TM\in\cont^1\di{\TTM\sp\TM}\,$
is called a \textit{spray}
if it is also a section of the bundle
$\,\TT\Tproj\in\cont^1\di{\TTM\sp\TM}\,$,
that is
if $\,\TT\Tproj\circ\BX=\Tproj_\TM\circ\BX=\ID\TM\,$.
\end{definition}
\begin{lemma}[Geodesics and sprays]\index{geodesic of a spray}
\label{lm: geospray}
Let $\,\spray\in\cont^1\di{\TM\sp\TTM}\,$ 
be a spray and $\,\Bv_\Bx\in\TANG_\Bx\MMM\,$ a tangent vector.
Then the base curve through $\,\Bx\in\MMM\,$
below the flow line of the spray through 
a vector $\,\Bv_\Bx\in\TANG_\Bx\MMM\,$
is a geodesic curve for any connection compatible with the spray,
i.e. such that $\,\Hlift_{\Bv_\Bx}\di{\Bv_\Bx}=\spray\di{\Bv_\Bx}\,$.
\end{lemma}
\proof
Let $\,\Bv_\ptau=\FL{\Hlift_\Bv}{\ptau,\pttt}\di{\Bv_\Bx}\,$
be the flow line of the spray through 
the vector $\,\Bv_\pttt=\Bv_\Bx\in\TANG_\Bx\MMM\,$.
The projected curve on the base manifold is then
$\,\cc_\ptau=(\TprojM\circ\FL{\spray}{\ptau,\pttt})\di{\Bv_\Bx}\,$,
with $\,\cc_\pttt=\Bx\in\MMM\,$.
Its velocity field 
$\,\Bv\in\cont^1\di{\II\sp\TM}\,$ is given by
$$\vcenter{\halign{
\hfil$#$&$#$\hfil&$#$\hfil\cr
\Bv_\ptau
&\,=\parder\xi\ptau\cc_\xi\di\Bx
=\parder\xi\ptau(\TprojM\circ\FL{\spray}{\xi,\pttt})\di{\Bv_\Bx}
\vspace{6pt}\cr
&\,=\TT\TprojM\punto\spray\di{\FL{\spray}{\ptau,\pttt}\di{\Bv_\Bx}}
=\Bpi_\TM\di{\spray\di{\FL{\spray}{\ptau,\pttt}\di{\Bv_\Bx}}}
=\FL{\spray}{\ptau,\pttt}\di{\Bv_\Bx}
\,,\cr}}$$
and $\,\spray\di{\Bv_\pttt}=\parder\ptau\pttt\parder\xi\ptau\cc_\xi\,$.
Being
$\,\Hlift_{\Bv_\Bx}\di{\Bv_\Bx}=\spray\di{\Bv_\Bx}\,$,
the formula for the time-covariant derivative yields:
$$\vcenter{\halign{
\hfil$#$&$#$\hfil&$#$\hfil\cr
\nablabar_{\pttt}\Bv
&\,=\parder\ptau\pttt\FL{\Hlift_\Bv}{\pttt,\ptau}\di{\Bv_\ptau}
=\parder \ptau\pttt\FL{\Hlift_\Bv}{\pttt,\ptau}\di{\FL{\spray}{\ptau,\pttt}\Bv_\pttt}
\vspace{6pt}\cr
&\,=\parder \ptau\pttt(\FL{\Hlift_\Bv}{\pttt,\ptau}\circ\FL{\spray}{\ptau,\pttt})\di{\Bv_\pttt}
=\spray\di{\Bv_\pttt}-\Hlift_{\Bv_\pttt}\di{\Bv_\pttt}=0
\,.\cr}}$$
Hence the curve $\,\cc\in\cont^1\di{\II\sp\MMM}\,$ is a geodesic.
\endprova
A similar proof shows that the base curve through
$\,\Bv_\Bx\in\TANG_\Bx\MMM\,$
below the tangent-flow line of a spray 
is the velocity field of a geodesic, 
in any connection compatible with the spray,
and that the velocity field of the base points of the line is a \Jacobi\ field
\cite{Michor1997}.

The next original result is the main contribution of this paper.
It provides, in the general context of fibre bundles,
the expression of the curvature in terms of covariant derivatives.

\begin{theorem}[Curvature and covariant derivatives]
\label{th: curvnice}
For a given section
$\,\Bs\in\cont^1\di{\MMM\sp\EEE}\,$
of a fibre bundle $\,\Bp\in\cont^1\di{\EEE\sp\MMM}\,$
and any pair of vector fields
$\,\Bu,\Bv\in\cont^1\di{\MMM\sp\TM}\,$,
the following identity holds on $\,\Bs\di\MMM\subset\EEE\,$:
$$\vcenter{\halign{
$#$\hfil&$#$\hfil&$#$\hfil\cr
\brack{\nablabar_\Bu}{\nablabar_\Bv}-\nablabar_{\brack\Bu\Bv}
+\brack{\Hlift_\Bu}{\Hlift_\Bv}-\Hlift_{\brack\Bu\Bv}
=\brack{\Hlift_\Bv}{\nablabar_\Bu}
+\brack{\nablabar_\Bv}{\Hlift_\Bu}=0
\,.\cr}}$$
Accordingly , the vertical-valued curvature two-form 
$\,\curvopbar_\Bx\di{\Bs}\di{\Bu,\Bv}\in\VERT_{\Bs\di\Bx}\EEE\,$
is given by
$$\curvopbar\di{\Bs}\di{\Bu,\Bv}=
\brack{\nablabar_\Bu}{\nablabar_\Bv}\di{\Bs}
-\nablabar_{\brack\Bu\Bv}\di{\Bs}\,.$$
\end{theorem}
\proof
By Lemma \ref{lm: relatedbrackets} we know that on $\,\Bs\di\MMM\subset\EEE\,$
it is
$\,\brack{\Tder\Bu}{\Tder\Bv}=\Tder{\brack\Bu\Bv}\,$.
By performing an extension of the natural derivatives,
e.g. by the foliation method envisaged in Lemma
\ref{lm: extfoliation}, the covariant derivatives of a section
$\,\Bs\in\cont^1\di{\MMM\sp\EEE}\,$ are consequently extended to 
(local) vector fields
$\,\nablabar_\Bu,\nablabar_\Bv\in\cont^1\di{\EEE\sp\VE}\,$. 
Then, being
$$\vcenter{\halign{
\hfil$#$&$#$\hfil&$#$\hfil\cr
\Tder\Bu=\nablabar_\Bu+\Hlift_\Bu\,,
\qquad
\Tder\Bv=\nablabar_\Bv+\Hlift_\Bv\,,
\qquad
\Tder{\brack\Bu\Bv}
=\nablabar_{\brack\Bu\Bv}+\Hlift_{\brack\Bu\Bv}
\,,\cr}}$$
by bilinearity of the \Lie\ bracket we get
$$\vcenter{\halign{
$#$\hfil&$#$\hfil&$#$\hfil\cr
\brack{\nablabar_\Bu+\Hlift_\Bu}{\nablabar_\Bv+\Hlift_\Bv}
&\,=\brack{\nablabar_\Bu}{\nablabar_\Bv}
+\brack{\Hlift_\Bu}{\Hlift_\Bv}
+\,\brack{\nablabar_\Bu}{\Hlift_\Bv}
+\brack{\Hlift_\Bu}{\nablabar_\Bv}
\vspace{8pt}\cr
&\,=\nablabar_{\brack\Bu\Bv}+\Hlift_{\brack\Bu\Bv}
\,,\cr}}$$
which, being 
$\,\brack{\Hlift_\Bu}{\Hlift_\Bv}-\Hlift_{\brack\Bu\Bv}
=\VPROJ\punto\brack{\Hlift_\Bu}{\Hlift_\Bv}\,$, can be written as:
$$\vcenter{\halign{
$#$\hfil&$#$\hfil&$#$\hfil\cr
\brack{\nablabar_\Bu}{\nablabar_\Bv}-\nablabar_{\brack\Bu\Bv}
+\VPROJ\punto\brack{\Hlift_\Bu}{\Hlift_\Bv}
=\brack{\Hlift_\Bv}{\nablabar_\Bu}+\brack{\nablabar_\Bv}{\Hlift_\Bu}
\,.\cr}}$$
The tensoriality of the curvature
$\,\VPROJ\punto\brack{\Hlift_\Bu}{\Hlift_\Bv}\,$,
as a function of the horizontal lifts
$\,\Hlift_\Bu\,$ and $\,\Hlift_\Bv\,$,
has the following implication.
Let the local vector fields
$\,\cF_\Bu^\Bx,\cF_\Bv^\Bx\in\cont^1\di{\EEE\sp\TE}\,$
be generated by dragging the vectors
$\,\Hlift_{\Bu_\Bx},\Hlift_{\Bv_\Bx}\in\TANG_{\Bs_\Bx}\EEE\,$ 
along the flows of the extended covariant derivatives
$\,\nablabar_\Bu,\nablabar_\Bv\in\cont^1\di{\EEE\sp\TE}\,$:
$$\vcenter{\halign{
$#$\hfil&$#$\hfil&$#$\hfil\cr
&\cF_\Bu^\Bx\circ\FL{\nablabar_\Bv}{\pttt}\equaldef
\TT\FL{\nablabar_\Bv}{\pttt}\circ\Hlift_{\Bu_\Bx}\,,
\vspace{8pt}\cr
&\cF_\Bv^\Bx\circ\FL{\nablabar_\Bu}{\pttt}\equaldef
\TT\FL{\nablabar_\Bu}{\pttt}\circ\Hlift_{\Bv_\Bx}
\,.\cr}}$$
By tensoriality, in evaluating the r.h.s. of the previous equality at a point
$\,\Bs\di\Bx\in\EEE\,$, the horizontal lifts
$\,\Hlift_\Bu,\Hlift_\Bv\in\cont^1\di{\EEE\sp\TE}\,$ 
can be substituted by the vector fields
$\,\cF_\Bu^\Bx,\cF_\Bv^\Bx\in\cont^1\di{\EEE\sp\TE}\,$.
Then, by definition:
$$\brack{\cF_\Bv^\Bx}{\nablabar_\Bu}_\Bx=0\,,\qquad 
\brack{\nablabar_\Bv}{\cF_\Bu^\Bx}_\Bx=0\,,$$ 
so that
$$\brack{\Hlift_\Bv}{\nablabar_\Bu}_\Bx+\brack{\nablabar_\Bv}{\Hlift_\Bu}_\Bx=
\brack{\cF_\Bv^\Bx}{\nablabar_\Bu}_\Bx+\brack{\nablabar_\Bv}{\cF_\Bu^\Bx}_\Bx=0\,.$$
The result holds for any extension of the natural derivatives
and the formula for the  curvature is independent of the extension,
since, by tensoriality, it depends only on the values of the covariant derivatives at
$\,\Bs\di\Bx\,$.
\endprova

\section{Connection on a vector bundle}
\label{sec: connvector}

Let us resume the peculiar properties of linear connections
on a vector bundle $\,\Bp\in\cont^1\di{\EEE\sp\MMM}\,$ to
infer the relevant special expression of the curvature form.
\begin{definition}[Linear connection]
\label{def: linhori}
In a vector bundle
$\,\Bp\in\cont^1\di{\EEE\sp\MMM}\,$
a connection is linear if the pair made of the horizontal lift 
$\,\Hlift_\Bv\in\cont^1\di{\EEE\sp\HE}\,$ and of the
vector field
$\,\Bv\in\cont^{1}\di{\MMM\sp\TM}\,$
is a linear vector bundle homomorphism 
from the vector bundle
$\,\Bp\in\cont^1\di{\EEE\sp\MMM}\,$
to the vector bundle
$\,\TT\Bp\in\cont^1\di{\TE\sp\TM}\,$.
This means that, given two sections
$\,\Bs_1,\Bs_2\in\cont^1\di{\MMM\sp\EEE}\,$, the following property of
$\,\Bp$-$\TT\Bp$-linearity holds:
$$\left\{\vcenter{\halign{
\hfil$#$&$#$\hfil&$#$\hfil\cr
&\,\Hlift_{\Bv_\Bx}\di{\Bs_1+_\Bp\Bs_2}=
\Hlift_{\Bv_\Bx}\Bs_1+_{\TT\Bp}\Hlift_{\Bv_\Bx}\Bs_2\,,
\vspace{8pt}\cr
&\,\Hlift_{\Bv_\Bx}\di{\aa\punto_\Bp\Bs}=
\aa\punto_{\TT\Bp}\Hlift_{\Bv_\Bx}\Bs\,,\perogni\aa\in\Re
\,.\cr}}\right.$$
\end{definition}
Being
$$\left\{\vcenter{\halign{
\hfil$#$&$#$\hfil&$#$\hfil\cr
&\,\TT_{\Bv_\Bx}\di{\Bs_1+_\Bp\Bs_2}=
\TT_{\Bv_\Bx}\Bs_1+_{\TT\Bp} \TT_{\Bv_\Bx}\Bs_2\,,
\vspace{8pt}\cr
&\, \TT_{\Bv_\Bx}\di{\aa\punto_\Bp\Bs}=
\aa\punto_{\TT\Bp} \TT_{\Bv_\Bx}\Bs\,,\perogni\aa\in\Re
\,,\cr}}\right.$$
the $\,\Bp$-$\TT\Bp$-linearity of the horizontal lift $\,\Hlift_{\Bv_\Bx}\,$
is equivalent to $\,\Bp$-$\TT\Bp$-linearity of the covariant derivative $\,\nablabar_{\Bv_\Bx}\,$:
$$\left\{\vcenter{\halign{
\hfil$#$&$#$\hfil&$#$\hfil\cr
&\,\nablabar_{\Bv_\Bx}\di{\Bs_1+_\Bp\Bs_2}=
\nablabar_{\Bv_\Bx}\Bs_1+_{\TT\Bp}\nablabar_{\Bv_\Bx}\Bs_2\,,
\vspace{8pt}\cr
&\,\nablabar_{\Bv_\Bx}\di{\aa\punto_\Bp\Bs}=
\aa\punto_{\TT\Bp}\nablabar_{\Bv_\Bx}\Bs\,,\perogni\aa\in\Re
\,.\cr}}\right.$$
The distinguishing feature with respect to a 
connection on a general fibre bundle
is that,
by the identification $\,\VE\simeq\EEE\,$,
the covariant derivative 
$\,\nablabar_\Bv\Bs\in\cont^1\di{\MMM\sp\VE}\,$
of a section
$\,\Bs\in\cont^1\di{\MMM\sp\EEE}\,$
along a vector field $\,\Bv\in\cont^1\di{\MMM\sp\TM}\,$
may be considered as a section
$\,\nabla_\Bv\Bs\in\cont^1\di{\MMM\sp\EEE}\,$
of the vector bundle
and the covariant derivative
$\,\nabla_\Bv\in\cont^1\di{\Bs\di\MMM\sp\EEE}\,$
as an operator.
The result stated below in proposition \ref{prop: Leibnizrulecov} 
makes appeal to this identification and is a basic property of the covariant derivative
in a linear connection (see e.g. \cite{Kobayashi1957}, \cite{KobayashiNomizu1963}).
\begin{proposition}[Leibniz rule for the covariant derivative]
\label{prop: Leibnizrulecov}
In a vector bundle
$\,\Bp\in\cont^1\di{\EEE\sp\MMM}\,$ 
endowed with a linear connection, 
the covariant derivative $\,\nabla_\Bv\in\cont^1\di{\Bs\di\MMM\sp\EEE}\,$
fulfils \Leibniz\ rule:
$$\nabla_{\Bv}\di{\ff\,\Bs}=(\nabla_{\Bv}\ff)\,\Bs+\ff\,(\nabla_{\Bv}\Bs)\,.$$
\end{proposition}
\proof
Let us recall from Lemma \ref{lm: covarLie} the expression
$\,\nablabar_\Bv\Bs
=\parder{\pttt}{0}
\FL{\Hlift_\Bv}{-\pttt}\circ\Bs\circ\FL{\Bv}{\pttt}\,$.
The linearity of the connection implies
that the flow $\,\FL{\Hlift_\Bv}{\pttt}\,$ is a one parameter family of automorphisms.
Then
$$\vcenter{\halign{
\hfil$#$&$#$\hfil&$#$\hfil&$#$\hfil&$#$\hfil\cr
\nablabar_{\Bv_\Bx}(\ff\Bs)
&\,=\parder{\pttt}{0}\ff\di{\FL\Bv\pttt\di\Bx}\punto
\FL{\Hlift_\Bv}{-\pttt}\di{\Bs\di{\FL\Bv\pttt\di\Bx}}
\,.\cr}}$$
To shorten the expressions we set
$\,\segno\ff_\pttt\equaldef\ff\di{\FL{\Bv}{\pttt}\di\Bx}\,$ and
$\,\segno\Bs_\pttt\equaldef\FL{\Hlift_\Bv}{-\pttt}\di{\Bs\di{\FL{\Bv}{\pttt}\di\Bx}}\,$
so that
$$\vcenter{\halign{
\hfil$#$&$\displaystyle#$\hfil&$#$\hfil&$#$\hfil&$#$\hfil\cr
\nablabar_{\Bv_\Bx}(\ff\Bs)
&\,=\parder\pttt0
\segno\ff_\pttt\punto\segno\Bs_\pttt
=\lim_{\pttt\to0}\inv\pttt
(\segno\ff_\pttt\segno\Bs_\pttt-\segno\ff_\pttt\segno\Bs_0
+\segno\ff_\pttt\segno\Bs_0-\segno\ff_0\segno\Bs_0)
\vspace{8pt}\cr
&\,=\lim_{\pttt\to0}\inv\pttt
(\segno\ff_\pttt(\segno\Bs_\pttt-\segno\Bs_0))
+\lim_{\pttt\to0}\inv\pttt
(\segno\ff_\pttt\segno\Bs_0-\segno\ff_0\segno\Bs_0)
\,.\cr}}$$
Hence, observing that
$$\vcenter{\halign{
\hfil$\displaystyle#$&$\displaystyle#$\hfil&$#$\hfil&$#$\hfil&$#$\hfil\cr
\lim_{\pttt\to0}\inv\pttt
(\segno\ff_\pttt(\segno\Bs_\pttt-\segno\Bs_0))
&\,=\segno\ff_0\,\parder\pttt0\segno\Bs_\pttt\,,
=\ff\di\Bx\,\nablabar_{\Bv_\Bx}\Bs\di\Bx\,,
\vspace{8pt}\cr
\lim_{\pttt\to0}\inv\pttt
(\segno\ff_\pttt\segno\Bs_0-\segno\ff_0\segno\Bs_0)
&\,=\lim_{\pttt\to0}\inv\pttt(\segno\ff_\pttt-\segno\ff_0)\,\segno\Bs_0
=(\nabla_{\Bv_\Bx}\ff)\,\Bs\di\Bx
\,,\cr}}$$
the result follows.
\endprova
In a vector bundle $\,\Bp\in\cont^1\di{\EEE\sp\MMM}\,$
the iterated and the second covariant derivatives
according to a given connection are meaningful. 
Hence, for any section
$\,\Bs\in\cont^1\di{\MMM\sp\EEE}\,$,
the curvature form may be written as
$$\vcenter{\halign{
\hfil$#$\hfil&$#$\hfil&$#$\hfil\cr
\curvoperator_{\Bs}\di{\Bu,\Bv}
&\,=(\nabla_\Bu\circ\nabla_\Bv-\nabla_\Bv\circ\nabla_\Bu
-\nabla_{\brack{\Bu}{\Bv}})\,\Bs
\vspace{6pt}\cr
&\,=(\nabla^2_{\Bu\Bv}-\nabla^2_{\Bv\Bu}
+\nabla_{\torstensor\di{\Bu,\Bv}})\,\Bs
\,,\cr}}$$
in terms of the second covariant derivative
$\,\nabla^2_{\Bu\Bv}\equaldef
\nabla_\Bu\circ\nabla_\Bv-\nabla_{\nabla_\Bu\Bv}\,$
and of the torsion
$\,\torstensor\di{\Bu,\Bv}\equaldef
\nabla_\Bu\Bv-\nabla_\Bv\Bu-\nabla_{\brack{\Bu}{\Bv}}\,$
which are both tensor fields.
Tensoriality may be proved by relying on
Leibniz rule to verify $\,\cont^\infty\di{\EEE\sp\Re}$-linearity.
\begin{remark}
We underline that on a fibre bundle, 
in writing  the formula:
$\,\curvoperator_{\Bs}\di{\Bu,\Bv}=
\brack{\nabla_\Bu}{\nabla_\Bv}\di{\Bs}
-\nabla_{\brack{\Bu}{\Bv}}\di{\Bs}\,$,
provided by Theorem \ref{th: curvnice}, the term
$\,\brack{\nabla_\Bu}{\nabla_\Bv}\di{\Bs}\,$ 
cannot be written as
$\,(\nabla_\Bu\circ\nabla_\Bv-\nabla_\Bv\circ\nabla_\Bu)\,\Bs\,$
since, being $\,\nabla_\Bu,\nabla_\Bv\in\cont^1\di{\EEE\sp\VE}\,$,
the compositions
$\,\nabla_\Bu\circ\nabla_\Bv\,$ and
$\,\nabla_\Bv\circ\nabla_\Bu\,$,
are not defined, unless the bundle is a vector bundle
and the identification 
$\,\VE\simeq\EEE\,$ can be made.
\end{remark}

\section{Conclusions}
\label{sec: Conclusions}

Connections on fibre bundles
and their torsion and curvature forms
are of primary importance in many
basic issues of mathematical physics, as witnessed by a 
vast number of contributions in literature (see e.g. \cite{MangiaSarda2000}).
The topic has been revisited here with the aim of providing a direct proof
to the relation between the integrability condition
provided by \Frobenius\ theorem and the 
expression of the curvature field in terms of covariant derivatives.
This result, in the general form provided here, appears to be new,
since classical treatments consider the special case of linear connections on vector bundles.
Our proof is based on the notion of natural derivative of a section, on
a suitable extension, by foliation, to a vector field in the tangent bundle, and on
a simple but powerful tensoriality argument.



\end{document}